\documentclass[12pt,draft]{amsart}
\usepackage{amsmath,amsthm,latexsym,amscd,amsbsy,amssymb}
 \relax

\chardef\bslash=`\\ 

\makeatletter
\def\verbatim{\interlinepenalty\@M \@verbatim
  \leftskip\@totalleftmargin\advance\leftskip2pc
  \frenchspacing\@vobeyspaces \@xverbatim}
\makeatother
\makeatletter
  \def\dgt@k{\dg@DX=-3 \dg@DY=2 \dg@SIZE=3} 
\makeatother

\makeatletter
  \def\dgt@kk{\dg@DX=3 \dg@DY=-1 \dg@SIZE=3}%
\makeatother

\theoremstyle{plain}
\newtheorem{thm}{Theorem}[section]
\newtheorem{cor}[thm]{Corollary}
\newtheorem{lem}[thm]{Lemma}
\newtheorem{pro}[thm]{Proposition}

\newtheorem*{A}{Theorem A}
\newtheorem*{B}{Theorem B}
\newtheorem*{C}{Theorem C}
\newtheorem*{D}{Theorem D}
\newtheorem*{E}{Theorem E}

\newtheorem*{F}{Theorem F}
\newtheorem*{G}{Theorem G}

\theoremstyle{definition}
\newtheorem{rem}[thm]{Remark}
\newtheorem{defin}[thm]{Definition}

\numberwithin{equation}{section}

\font\f=msbm10

\begin{document}


\title[Topological $\text{AE}(0)$-groups]{Topological $\text{AE}(0)$-groups}
\author{Alex Chigogidze}
\address{Department of Mathematics and Statistics,
University of Saskatche\-wan,
McLean Hall, 106 Wiggins Road, Saskatoon, SK, S7N 5E6,
Canada}
\email{chigogid@math.usask.ca}
\thanks{Author was partially supported by NSERC research grant.}



\begin{abstract}{We investigate topological
$\text{AE}(0)$-groups class of which
contains the class of Polish groups as well as the class of
all locally compact groups. We establish the existence
of an universal $\text{AE}(0)$-group of a given weight
as well as the existence of an universal action of
$\text{AE}(0)$-group of a given weight on a $\text{AE}(0)$-space
of the same weight. A complete characterization of closed
subgroups of powers of the symmetric group $S_{\infty}$ is obtained.
It is also shown that every $\text{AE}(0)$-group is Baire
isomorphic to the product of Polish groups. These results
are obtained by using the spectral descriptions of
$\text{AE}(0)$-groups which are presented in Section \ref{S:spectral}.}
\end{abstract}

\maketitle
\markboth{A.~Chigogidze}{Topological $\text{AE}(0)$-groups}

\section{Introduction}\label{S:intro}

One of the main structure theorems for compact groups
(see \cite[ Chapters 6 \& 9]{hofmor98}) can be formulated as follows.

\begin{A}[\cite{hofmor98}, Theorem 9.24(ii)]\label{T:A}
Let $G$ be a connected compact group. Then there exists
a continuous homomorphism 
\[ p \colon Z_{0}(G) \times \prod\{ L_{t} \colon t \in T\} \to G ,\]
where $Z_{0}(G)$ stands for the identity component of the
center of $G$, and $L_{t}$ is a simple, connected and simply
connected compact Lie group, $t \in T$, such that $\ker (p)$
is a zero-dimensional central subgroup of
$\displaystyle Z_{0}(G) \times \prod\{ L_{t} \colon t \in T\}$.
\end{A}


The above statement clearly shows that the classes of
zero-dimensional groups,
abelian groups and simple, simply connected Lie groups
play a central role in the
general theory of compact groups. Recently these classes of
topological groups
have been studied from the point of view 
of absolute extensors in dimension $n$ (see \cite{chibook96}
for a comprehensive introduction into the theory of absolute
extensors in dimension $n$). The following three statements
show that such an approach is quite effective.

\begin{B}[\cite{bellchi93}]\label{T:B}
The following conditions are equivalent for a
zero-dimen\-si\-o\-nal topological group $G$:
\begin{itemize}
\item[(a)]
$G$ is topologically equivalent to the product
$(\text{\f Z}_{2})^{\tau} \times \text{\f Z}^{\kappa}$.
\item[(b)]
$G$ is an $\text{AE}(0)$-space.
\end{itemize}
\end{B}

\begin{C}[\cite{chi993}, Theorem E]\label{T:C}
The following conditions are equivalent for a compact
abelian group $G$:
\begin{itemize}
\item[(a)]
$G$ is a torus group (both in topological and algebraic senses).
\item[(b)]
$G$ is an $\text{AE}(1)$-compactum.
\end{itemize}
\end{C}

\begin{D}[\cite{chi972}, Corollary 1]\label{T:D}
The following conditions are equivalent for a non-trivial
compact group $G$:
\begin{itemize}
\item[(a)]
$G$ is a simple, connected and simply connected Lie group.
\item[(b)]
$G$ is an $\text{AE}(2)$-group with $\pi_{3}(G) = \text{\f Z}$.
\end{itemize}
\end{D}

The full classification problem for 
non-abelian (see Theorem C) and for non-simply connected
compact $\text{AE}(1)$-groups remains open. On the other hand,
the following two statements provide a complete classification of
simply connected compact $\text{AE}(1)$-groups.

\begin{E}[\cite{chi972}, Theorem C]\label{T:E}
The following conditions are equivalent for a compact group $G$:
\begin{itemize}
\item[(a)]
$G$ is a simply connected $\text{AE}(1)$-compactum.
\item[(b)]
$G$ is an $\text{AE}(2)$-compactum.
\item[(c)]
$G$ is an $\text{AE}(3)$-compactum.
\item[(d)]
$G$ is a product of simple, connected and simply connected
compact Lie groups.
\end{itemize}
\end{E}

\begin{F}[\cite{chi972}, Corollary 2]\label{T:F}
There is no non-trivial compact
$\text{AE}(4)$-group.
\end{F}

We complete this brief survey by pointing out that, as shows the
following statement,
for locally compact 
groups the restriction of being $\text{AE}(0)$-group is purely formal.

\begin{G}[Pontryagin-Haydon, \cite{pont46},
\cite{hay74}; \cite{chibook96}]\label{T:G}
Every locally compact group is an $\text{AE}(0)$-space.
\end{G}

Below we study $\text{AE}(0)$-groups. 
The class of $\text{AE}(0)$-groups contains the class of all
(generally speaking, non-metrizable)
locally compact groups (Theorem G) as well as the class of
all Polish groups. Actually the class of Polish groups coincides
with the class of metrizable $\text{AE}(0)$-groups and forms a
foundation of the entire theory of $\text{AE}(0)$-groups. We 
hope that results presented in Section \ref{S:spectral} do
indicate a potential for a non-trivial
theory of $\text{AE}(0)$-groups which unifies and generalizes
theories of locally compact and Polish groups (thus offering a
possible approach to the
corresponding question posed in \cite{pest99}).

In Section \ref{S:spectral} we present a spectral characterization of
$\text{AE}(0)$-groups in terms of well ordered continuous
inverse spectra (Theorem \ref{T:ae0}). This characterization
states that a non-metrizable topological group $G$ of weight $\tau$
is a $\text{AE}(0)$-group if and
only if it is the limit of a well ordered continuous inverse system
${\mathcal S}_{G} = \{ G_{\alpha}, p_{\alpha}^{\alpha +1},
\alpha < \tau \}$ of length $\tau$, consisting of
$\text{AE}(0)$-groups $G_{\alpha}$ and $0$-soft
limit homomorphisms 
$p_{\alpha}^{\alpha +1} \colon G_{\alpha +1} \to G_{\alpha}$,
$\alpha < \tau$, so that $G_{0}$ is a Polish group and each homomorphism
$p_{\alpha}^{\alpha +1}$, $\alpha < \tau$, has a Polish kernel.

Obviously this result can not be accepted as the one providing a
satisfactory reduction of the non-metrizable case to the Polish
one. Of course, everything is fine if the weight of $G$ is
$\omega_{1}$ -- in such a case all $G_{\alpha}$'s,
$\alpha < \omega_{1}$, (and not only the very first one,
i.e. $G_{0}$) are indeed Polish. But if the weight of $G$
is greater than $\omega_{1}$, then all $G_{\alpha}$'s, with
$\alpha \geq \omega_{1}$, are non-metrizable. 

In order to achieve our final goal and complete the reduction, we
analyze $0$-soft homomorphisms with Polish kernels between
(generally speaking, non-metrizable)
$\text{AE}(0)$-groups. A characterization of such
homomorphisms, which is recorded in Proposition \ref{P:polishkernel},
states that a $0$-soft homomorphism
$f \colon G \to L$ of $\text{AE}(0)$-groups has a Polish kernel
if and only if there exists a pullback diagram
\[
\begin{CD}
G @>f>> L\\
@V{p}VV @VV{q}V\\
G_{0} @>f_{0}>> L_{0},
\end{CD}
\]

\noindent where $G_{0}$ and $L_{0}$ are Polish groups
and the homomorphisms $p \colon G \to G_{0}$ and
$q \colon L \to L_{0}$ are $0$-soft. Theorem
\ref{T:ae0} and Proposition \ref{P:polishkernel}
together complete the required reduction. 

Polish groups and their actions have been
extensively studied in a variety of directions (ergodic theory,
group representations,
operator algebras; see \cite{bekech96} for further
discussion and references). Some of the central themes of
the theory of Polish groups counterparts of which (for
arbitrary $\text{AE}(0)$-groups) are considered below  are the
existence of universal groups, the existence of universal actions and
characterization of closed subgroups of the symmetric group $S_{\infty}$.

In Section \ref{S:applications} we use above mentioned spectral
characterizations and present extensions of some of these results
for $\text{AE}(0)$-groups. We prove the existence of  universal
$\text{AE}(0)$-groups of a given weight (Proposition \ref{P:univ})
and the existence of universal actions of
$\text{AE}(0)$-groups of a given weight 
on compact $\text{AE}(0)$-spaces of the same
weight (Theorem \ref{T:actionuniv}).  

Theorem \ref{T:closed} characterizes $\text{AE}(0)$-groups which
are isomorphic to closed subgroups of powers $S_{\infty}^{\tau}$ of 
the symmetric group $S_{\infty}$ - the group
of all bijections of ${\mathbb N}$ under the relative topology
inherited from ${\mathbb N}^{\mathbb N}$. This result extends the
corresponding observation \cite[Theorem 1.5.1]{bekech96} for the
group $S_{\infty}$ itself. As a corollary (Corollary
\ref{C:polish}) we note that if a
Polish group $G$ can be embedded (as a closed subgroup) into
$S_{\infty}^{\tau}$ for some $\tau$, then $G$ can be embedded
into $S_{\infty}$ as well. In light of
\cite{dough94} this shows that there exist zero-dimensional
Polish groups
which can not be embedded into $S_{\infty}^{\tau}$ as closed
subgroups for any cardinal number $\tau$.

Finally we use Theorem \ref{T:ae0} to prove (Theorem \ref{T:baire})
that every $\text{AE}(0)$-group is Baire isomorphic to the product of Polish groups.
In light of Theorem G this result appears to be new even for compact groups.


\section{Preliminaries}\label{S:pre}
All topological spaces below are assumed to be Tychonov (i.e.
completely regular and Hausdorff) and all maps (except
in Subsection \ref{SS:baire}) are continuous. We consider
only Lebesgue dimension $\dim$. Definitions of concepts related to inverse spectra
can be found in \cite{chibook96}. ${\mathbb R}$
denotes the real line and ${\mathbb Q}$ stand for the Hilbert cube.

\subsection{Definitions of $\text{AE}(n)$-spaces and $n$-soft maps}\label{SS:aen}
A comprehensive introduction into general theory
of $\text{AE}(n)$-spaces and $n$-soft maps can be found in \cite{chibook96}.
$C(X)$ denotes the set of all continous real-valued functions defined on $X$.

\begin{defin}\label{D:6.1.3}
A space $X$ is called an absolute  extensor in dimension $n$
(shortly, $\text{AE}(n)$-space), $n=0,1,\dots$, if for each at
most $n$-dimensional space $Z$ and each subspace $Z_0$ of $Z$,
any map $f \colon Z_0 \to X$, such that
$C(f)(C(X)) \subseteq \{ \varphi |Z_{0} \colon \varphi \in C(Z) \}$,
can be extended to $Z$.
\end{defin}

For compact spaces this definition is equivalent to the standard one.

\begin{pro}\label{P:6.1.11}
A compact space $X$ is a $\text{AE}(n)$-space
if and only if for each at most $n$-di\-men\-si\-o\-nal
compactum $Z$ and for each closed subspace $Z_0$ of $Z$,
any map $f \colon Z_0 \to X$ has an extension to  $Z$.
\end{pro}

It is known \cite[Chapter 6]{chibook96} that the
class of metrizable $\text{AE}(0)$-spaces
coincides with the class of Polish spaces. Every
$\text{AE}(0)$-spaces has a countable Suslin number. 

\begin{defin}\label{D:6.1.16}
A map $f \colon X \to Y$ between $\text{AE}(n)$-spaces
is $n$-soft if and only if for each at most $n$-dimensional
realcompact space $Z$, for its closed subspace $Z_0$, and for
any two maps $g \colon Z_0 \to X$ and $h \colon Z \to Y$ such that
$f\circ g = h|Z_0$ and
$C(g)(C(X)) \subseteq \{ \varphi |Z_{0} \colon \varphi \in C(Z)\}$,
there exists a map $k \colon Z \to X$ such that
$k|Z_0 = g$ and $f\circ k = h$.
\end{defin}

It is easy to check that $X$ is $\text{AE}(n)$-space if and only if
the constant map $X \to \{ \operatorname{pt}\}$ is $n$-soft.
It is important to note that every $0$-soft map between $\text{AE}(0)$-spaces
is surjective and open (\cite[Lemma 6.1.13 \& Proposition 6.1.26]{chibook96}) and that
for surjections between Polish spaces the converse of this fact is also true.

We say (see \cite[Section 6.3]{chibook96}) that
a map $f \colon X \to Y$ has a Polish kernel if there exists a Polish
space $P$ such that
$X$ is $C$-embedded in the product $Y \times P$ so that $f$ coincides
with the restriction $\pi_{Y}|X$ of the projection
$\pi_{Y} \colon Y \times P \to Y$. Obviously any map between
Polish spaces has a Polish kernel.

\subsection{Set-theoretical preliminaries}\label{SS:set}
Let $A$  be a partially ordered directed set (i.e.
for every two elements  $\alpha ,\beta \in A$  there exists
an element  $\gamma \in A$  such that  $\gamma \geq \alpha$ 
and  $\gamma \geq \beta$). We say that a subset
$A_1 \subseteq A$ of $A$ majorates another subset
$A_2 \subseteq A$ of $A$ if for each element $\alpha_2 \in A_2$
there exists an element $\alpha_1 \in A_1$ such that
$\alpha_1 \geq \alpha_2$. A subset which majorates $A$
is called cofinal in $A$. A subset of  $A$  is said to
be a chain if every two elements of it are comparable.
The symbol $\sup B$ , where  $B \subseteq A$, denotes the
lower upper bound of $B$ (if such an element exists in $A$).
Let now $\tau$ be an infinite cardinal number. A subset $B$
of $A$  is said to be $\tau$-closed in $A$ if for each chain
$C \subseteq B$, with ${\mid}C{\mid} \leq \tau$, we have
$\sup C \in B$, whenever the element $\sup C$ exists in $A$.
Finally, a directed set $A$ is said to be $\tau$-complete
if for each chain $B$ of elements of $A$ with
${\mid}C{\mid} \leq \tau$, there exists an element
$\sup C$ in $A$. 

The standard example of a $\tau$-complete set can be obtained
as follows. For an arbitrary set $A$ let $\exp A$ denote, as usual,
the collection of all subsets of $A$. There is a natural partial
order on $\exp A$: $A_1 \geq A_2$ if and only if $A_1 \supseteq A_2$.
With this partial order $\exp A$ becomes a directed set.
If we consider only those subsets of the set $A$ which have
cardinality $\leq \tau$, then the corresponding subcollection
of $\exp A$, denoted by $\exp_{\tau}A$, serves as a basic
example of a $\tau$-complete set. Proofs of the following
statements can be found in \cite{chibook96}.

\begin{pro}\label{P:3.1.1}
Let  $\{ A_{t} : t \in T \}$ be a collection of $\tau$-closed and
cofinal subsets of a $\tau$-complete set $A$. If
$\mid T\mid \leq \tau$, then the intersection
$\cap \{ A_{t}: t \in T \}$ is also cofinal
(in particular, non-empty) and $\tau$-closed in $A$ .
\end{pro}

\begin{cor}\label{C:3.1.2}
For each subset $B$, with  $\mid B \mid \leq \tau$, of a
$\tau$-complete set $A$ there exists an element $\gamma \in A$
such that  $\gamma \geq \beta$  for each  $\beta \in B$ .
\end{cor}

\begin{pro}\label{P:search}
Let  $A$  be a $\tau$-complete set, 
$L \subseteq A^2$, and suppose the following three
conditions are satisfied:
\begin{description}
\item[Existence] For each $\alpha \in A$ there exists
$\beta \in A$  such that  $(\alpha ,\beta ) \in L$.
\item[Majorantness] If  $(\alpha ,\beta ) \in L$  and
$\gamma \geq \beta$, then  $(\alpha ,\gamma ) \in L$.
\item[$\tau$-closeness] Let $\{ \alpha_{t} : t \in T \}$
be a chain in $A$ with $|T| \leq \tau$. If
$(\alpha_{t}, \beta ) \in L$ for some
$\beta \in A$ and each $t \in T$, then
$(\alpha ,\beta ) \in L$ where $\alpha =
\sup \{\alpha_{t} \colon t \in T \}$.
\end{description}
   Then the set of all  $L$-{\em reflexive} elements of 
$A$ (an element $\alpha \in A$ is $L$-reflexive if
$(\alpha ,\alpha ) \in L$)  is cofinal and $\tau$-closed in $A$.
\end{pro}

\subsection{Baire sets and Baire isomorphisms}\label{SS:bairesets}
Recall that elements of the $\sigma$-algebra generated
by functionally open subsets of a space $X$ are called Baire sets of $X$.
A map $f \colon X \to Y$ is a Baire map if inverse images of Baire
sets are Baire sets. A bijection $f \colon X \to Y$ is Baire
isomorphism if both $f$ and $f^{-1}$ are Baire maps.

The following statement (\cite[Proposition 2.5]{chi994}) is
used in the proof of Theorem \ref{T:baire}.

\begin{pro}\label{P:bb}
Let $f \colon X \to Y$ be a $0$-soft map of
$\text{AE}(0)$-spaces. Then there exists a Baire map
$g \colon Y \to X$ such that $f \circ g = \operatorname{id}_{Y}$.
\end{pro}

If $f$ is a continuous homomorphism of Polish groups, then
the existence of such a $g$ has been observed by Dixmier
\cite[Theorem 1.2.4]{bekech96}, \cite[Theorem 12.17]{kech95}).


\section{$\text{AE}(0)$-groups and actions of
$\text{AE}(0)$-groups -- spectral representations}\label{S:spectral}

In this section we present spectral characterizations of
$\text{AE}(0)$-groups. We also present a spectral description
of actions of $\text{AE}(0)$-groups on $\text{AE}(0)$-spaces. 

\begin{pro}\label{P:B8.2.1}
Each  topological $\text{AE}(0)$-group $G$ is topologically
and algebraically isomorphic to the limit of a factorizing
$\omega$-spectrum ${\mathcal S}_{G} = \{ G_{\alpha},
p_{\alpha}^{\beta}, A\}$
consisting of Polish groups $G_{\alpha}$, $\alpha \in A$, and $0$-soft
limit homomorphisms $p_{\alpha} \colon G \to G_{\alpha}$, $\alpha \in A$.
In particular, all short projections
$p_{\alpha}^{\beta} \colon G_{\beta} \to G_{\alpha}$,
$\alpha \leq \beta$, $\alpha ,\beta \in A$, are $0$-soft homomorphisms.
\end{pro}
\begin{proof}
By \cite[Theorem 6.3.2 or Proposition 6.3.5]{chibook96},
the space $G$ can be
represented as the limit space of a factorizing $\omega$-spectrum
${\mathcal S} = \{ G_{\alpha}, p_{\alpha}^{\beta}, \widetilde{A}\}$
consisting
of Polish spaces (i.e. $\text{AE}(0)$-spaces of countable weight)
and $0$-soft limit
projections. Let us show that this spectrum contains
$\omega$-closed and cofinal subspectrum consisting of
topological groups and limit projections that are
(continuous) homomorphisms. 

Let $\mu \colon G \times G \to G$ and $\nu \colon G \to G$ be
continuous operations of multiplication and inversion given on
$G$ as on a topological group. We apply \cite[Theorem 1.3.6]{chibook96} to
both $\mu$ and $\nu$. First consider the multiplication.
Clearly, $G \times G$ is the limit space
of the spectrum\\
\[ {\mathcal S} \times {\mathcal S} =
\{ G_{\alpha} \times G_{\alpha}, p_{\alpha}^{\beta}
\times p_{\alpha}^{\beta}, \widetilde{A} \} .\]

\noindent All projections of the spectrum
${\mathcal S} \times {\mathcal S}$ are $0$-soft, and hence,
by \cite[Proposition 6.1.26]{chibook96}, open. The Suslin number
of the product $G \times G$ is obviously countable
(see \cite[proposition 6.1.8]{chibook96}). Consequently,
by \cite[Proposition 1.3.3]{chibook96}, the spectrum
${\mathcal S} \times {\mathcal S}$ is factorizing. Next
we apply \cite[Theorem 1.3.6]{chibook96} to the spectra
${\mathcal S} \times {\mathcal S}$, $\mathcal S$ and to the
map $\mu$ between their limit spaces. Then we get a
$\omega$-closed and cofinal subset $A_{\mu}$ of $\widetilde{A}$
such that for
each $\alpha \in A_{\mu}$ there exists a continuous map
$\mu _{\alpha} \colon G_{\alpha} \times G_{\alpha} \to G_{\alpha}$
such that the diagram\\

\[
\begin{CD}
G \times G @>\mu>> G\\
@V{p_{\alpha}\times p_{\alpha}}VV @VV{p_{\alpha}}V\\
G_{\alpha}\times G_{\alpha} @>\mu_{\alpha}>> G_{\alpha}
\end{CD}
\]
\bigskip

\noindent commutes. In other words,
$p_{\alpha}\circ \mu = \mu _{\alpha}\circ (p_{\alpha} \times p_{\alpha})$
for each $\alpha \in A_{\mu}$.

Next consider a continuous inversion $\nu \colon G \to G$. By
\cite[Theorem 1.3.6]{chibook96}, applied to $\nu$ and the
spectrum ${\mathcal S}$,
there exists a $\omega$-closed and cofinal subset
$A_{\nu}$ of $\widetilde{A}$ such that for
each $\alpha \in A_{\nu}$ there exists a continuous map
$\nu _{\alpha} \colon G_{\alpha}  \to G_{\alpha}$
such that the diagram\\

\[
\begin{CD}
G  @>\nu>> G\\
@V{p_{\alpha}}VV @VV{p_{\alpha}}V\\
G_{\alpha} @>\nu_{\alpha}>> G_{\alpha}
\end{CD}
\]
\bigskip

\noindent commutes. In other words,
$p_{\alpha}\circ \nu = \nu_{\alpha} \circ p_{\alpha}$
for each $\alpha \in A_{\nu}$.

By Proposition \ref{P:3.1.1}, the intersection
$A = A_{\mu} \cap A_{\nu}$ is still $\omega$-closed and cofinal
in $\widetilde{A}$. This guarantees that $G$ is topologically
and algebraically isomorphic to the limit of the factorizing
$\omega$-spectrum
${\mathcal S}_{G} = \{ G_{\alpha}, p_{\alpha}^{\beta} , A\}$.
Also for each $\alpha \in A$ we have two maps
\[ \mu_{\alpha} \colon G_{\alpha} \times G_{\alpha} \to
G_{\alpha}\;\;\text{and}\;\; \nu_{\alpha} \colon G_{\alpha}
\to G_{\alpha} \]

\noindent which allow us to define

\begin{itemize}
\item[(a)]
a continuous multiplication operation on 
$G_{\alpha}$ by letting\\
\[ x_{\alpha}\cdot y_{\alpha} = \mu _{\alpha}(x_{\alpha},
y_{\alpha})\;\;
\text{for each}\;\; (x_{\alpha},y_{\alpha}) \in G_{\alpha}
\times G_{\alpha};\]
\medskip
\item[(b)]
a continuous inversion on $G_{\alpha}$ by letting 

\[ x_{\alpha}^{-1} = \nu_{\alpha}(x_{\alpha}) \;\;
\text{for each}\;\; x_{\alpha} \in G_{\alpha} .\]
\medskip
\end{itemize}

It is easy to see that $G_{\alpha},\; \alpha \in A$, becomes a
topological group
with respect to these operations. Moreover, for each $\alpha \in A$
the limit projection $p_{\alpha} \colon G \to G_{\alpha}$ becomes a
homomorphism with respect to the above defined operations.
\end{proof} 

\begin{cor}\label{C:dim}
Let $G$ be a $\text{AE}(0)$-group such that $\dim G \leq n$. Then
$G$ is topologically
and algebraically isomorphic to the limit of a factorizing
$\omega$-spectrum ${\mathcal S}_{G} = \{ G_{\alpha},
p_{\alpha}^{\beta}, A\}$
consisting of at most $n$-dimensional Polish groups $G_{\alpha}$,
$\alpha \in A$, and $0$-soft
limit homomorphisms $p_{\alpha} \colon G \to G_{\alpha}$, $\alpha \in A$.
\end{cor}
\begin{proof}
By Proposition \ref{P:B8.2.1}, $G = \lim {\mathcal S}_{1}$,
where ${\mathcal S}_{1} = \{ G_{\alpha}, p_{\alpha}^{\beta}, A_{1}\}$
is a factorizing $\omega$-spectrum consisting of Polish groups
and $0$-soft limit homomorphisms. By \cite[Theorem 1.3.10]{chibook96},
$G = \lim {\mathcal S}_{2}$, where
${\mathcal S}_{2} = \{ G_{\alpha}, p_{\alpha}^{\beta}, A_{2}\}$
is a factorizing $\omega$-spectrum consisting of at most
$n$-dimensional Polish spaces. It follows directly from the
proofs of Proposition \ref{P:B8.2.1} and
\cite[Theorem 1.3.10]{chibook96}, that both indexing sets
$A_{1}$ and $A_{2}$ can be assumed to be closed and
$\omega$-complete subsets of a $\omega$-complete set $B$.
By Proposition
\ref{P:3.1.1}, $A = A_{1} \cap A_{2}$ is still cofinal and
$\omega$-closed subset of $B$.
Consequesntly, the factorizing $\omega$-spectrum
${\mathcal S}_{G} = \{ G_{\alpha}, p_{\alpha}^{\beta}, A\}$
consists of at most $n$-dimensional Polish groups and $0$-soft
limit homomorphisms. It only remains to note that
$G = \lim{\mathcal S}_{G}$.
\end{proof}

\begin{cor}\label{C:emb}
Let $\tau \geq \omega$. Every $\text{AE}(0)$-group of
weight $\tau \geq \omega$
is topologically and algebraically isomorphic to a closed
and $\text{C}$-embedded
subgroup of the product
$\displaystyle \prod\{ G_{t} \colon t \in T\}$, where $G_{t}$,
$t \in T$, is a Polish group and $|T| = \tau$.
\end{cor}
\begin{proof}
Let $G$ be an $\text{AE}(0)$-group of weight $\tau$.
If $\tau = \omega$, then
$G$ is Polish (see \cite[Chapter 6]{chibook96} and consequently
there is nothing to prove.

If $\tau > \omega$, then by Proposition \ref{P:B8.2.1}, $G$ is
topologically and algebraically 
isomorphic to the limit of a factorizing $\omega$-spectrum
${\mathcal S}_{G} = \{ G_{t}, p_{t}^{t^{\prime}}, T\}$
with $|T|  = \tau$. 
Clearly $\lim {\mathcal S}_{G}$ is isomorphic to a closed subgroup
of the product
$\displaystyle \prod\{ G_{t}, t \in A\}$. Since the spectrum
${\mathcal S}_{G}$ is factorizing,
it follows that $\lim {\mathcal S}_{G}$ is $\text{C}$-embedded
in $\displaystyle \prod\{ G_{t}, t \in T\}$.
\end{proof}


\begin{thm}\label{T:ae0}
Let $G$ be a topological group of weight $\tau > \omega$.
Then the following conditions are equivalent:
\begin{itemize}
\item[(a)]
$G$ is a $\text{AE}(0)$-group.
\item[(b)]
There exists a well-ordered inverse spectrum
${\mathcal S}_{G} = \{ G_{\alpha},
p_{\alpha}^{\alpha +1},\tau \}$ satisfying the following properties:
\begin{enumerate}
\item
$G_{\alpha}$ is a $\text{AE}(0)$-group and $p^{\alpha +1}_{\alpha}
\colon G_{\alpha +1} \to G_{\alpha}$ is a $0$-soft
homomorphism with Polish kernel, $\alpha < \tau$.
\item
If $\beta < \tau$ is a limit ordinal, then the diagonal product 
\[\triangle\{ p_{\alpha}^{\beta} \colon \alpha < \beta\}
\colon G_{\beta} \to \lim\{ G_{\alpha},
p_{\alpha}^{\alpha+1}, \alpha < \beta \}\]
is a topological and algebraic isomorphism.
\item
$G$ is topologically and algebraically isomorphic to
$\lim {\mathcal S}_{G}$.
\item
$G_{0}$ is a Polish group. 
\end{enumerate}
\end{itemize}
\end{thm}
\begin{proof}
(a) $\Longrightarrow$ (b). By Corollary \ref{C:emb}, we may
assume that $G$ is a closed and $\text{C}$-embedded subgroup
of the product $\displaystyle \prod\{ X_{a} \colon a \in A\}$,
$|A| = \tau$, of Polish groups $X_{a}$, $a \in A$. There exists
a proper, functionally closed and $0$-invertible map
$f \colon Y \to \prod\{ X_{a} \colon a \in A \}$, where $Y$
is a spectrally complete (see \cite[p.247]{chibook96})
realcompact space of weight $\tau$ and dimension $\dim Y = 0$
(see \cite[Proposition 6.2.13]{chibook96} for details).
Consider the inverse image $f^{-1}(G) \subseteq Y$ of $G$ and
the map $f|f^{-1} \colon f^{-1}(G) \to G$. Since $G$ is
$\text{C}$-embedded in the product $\displaystyle \prod\{
X_{a} \colon a \in A\}$, since $\dim Y = 0$ and since $G$,
according to (a), is an $AE(0)$-space, there exists a map
$g \colon Y \to G$ such that $g|f^{-1}(G) = f|f^{-1}(G)$.

Next let us denote by 
\[ \pi_{B} \colon \prod\{ X_{a} \colon a \in A\} \to
\prod\{ X_{a} \colon a \in B\}\] 
and
\[\pi^{B}_{C} \colon \prod\{ X_{a} \colon a \in B\} \to
\prod\{ X_{a} \colon a \in C \}\]
the natural projections onto the corresponding subproducts
($C \subseteq B \subseteq A$).
We call a subset $B \subseteq A$ admissible (compare
with the proof of \cite[Theorem 6.3.1]{chibook96}) if the
following equality
\[ \pi_{B}(g(f^{-1}(x))) = \pi_{B}(x)\]
is true for each point $x \in \pi_{B}^{-1}\left(\pi_{B}(G)
\right)$. We need the following properties of admissible sets.

{\em Claim 1. The union of arbitrary collection of admissible
sets is admissible}. 

Indeed let $\{ B_{t} \colon t \in T\}$ be a collection of
admissible sets and $B = \cup \{ B_{t} \colon t \in T\}$.
Let $x \in \pi_{B}^{-1}\left(\pi_{B}(G)\right)$. Clearly
$x \in \pi_{B_{t}}^{-1}\left(\pi_{B_{t}}(G)\right)$ for each
$t \in T$ and consequently 
\[ \pi_{B_{t}}(g(f^{-1}(x))) = \pi_{B_{t}}(x) \;\;\text{for
each}\;\; t \in T .\]
Obviously, $\pi_{B}(x) \in \pi_{B}(g(f^{-1}(x)))$ and it
therefore suffices to show that the set $\pi_{B}(g(f^{-1}(x)))$
contains only one point. Assuming that there is a point
$y \in \pi_{B}(g(f^{-1}(x)))$ such that $y \neq \pi_{B}(x)$
we conclude (having in mind that $B = \cup\{ B_{t} \colon
t \in T\}$) that there must be an index $t \in T$ such that
$\pi_{B_{t}}^{B}(y) \neq \pi_{B_{t}}^{B}\left(\pi_{B}(x)\right)$.
But this is impossible
\[ \pi_{B_{t}}^{B}(y) \in \pi_{B_{t}}^{B}\left(\pi_{B}
(g(f^{-1}(x)))\right) = \pi_{B_{t}}(g(f^{-1}(x))) = \pi_{B_{t}}(x) = \pi_{B_{t}}^{B}\left(\pi_{B}(x)\right) .\]

{\em Claim 2. If $B \subseteq A$ is admissible, then the
restriction $\pi_{B}|G \colon G \to \pi_{B}(G)$ is $0$-soft}.

Let $\varphi \colon Z \to \pi_{B}(G)$ and
$\varphi_{0} \colon Z_{0} \to G$ be two maps defined on a
realcompact space $Z$, with $\dim Z = 0$, and its closed
subset $Z_{0}$ respectively. Assume that
$\pi_{B}\varphi_{0} = \varphi |Z_{0}$ and
$C(\varphi_{0})(C(G)) \subseteq C(Z)|Z_{0}$. We wish to
construct a map $\phi \colon Z \to G$ such that
$\phi |Z_{0} = \varphi_{0}$ and $\pi_{B}\phi = \varphi$, i.e.
$\phi$ makes the diagram

\begin{picture}(300,150)
\put(80,110){$Z$}
\put(80,10){$Z_{0}$}
\put(80,100){$\cap$}
\put(87,102){\vector(0,-1){79}}
\put(65,60){$\operatorname{incl}$}
\put(95,17){\vector(1,1){90}}
\put(95,112){\vector(1,0){90}}
\put(95,12){\vector(1,0){80}}
\put(190,110){$G$}
\put(180,10){$\pi_{B}(G)$}
\put(195,107){\vector(0,-1){83}}
\put(135,118){$\varphi_{0}$}
\put(135,17){$\varphi$}
\put(135,67){$\phi$}
\put(198,60){$\pi_{B}|G$}
\end{picture}

\noindent commutative. Since, according to our choice,
all $X_{a}$'s are $\text{AE}(0)$-spaces (recall that each
$X_{a}$ is a Polish space), so is the product
$\prod\{ X_{a} \colon a \in A-B\}$. This implies the $0$-softness
of the projection $\pi_{B}$ and hence of its restriction 
\[ \pi_{B}|\pi_{B}^{-1}\left(\pi_{B}(G)\right) \colon
\pi_{B}^{-1}\left(\pi_{B}(G)\right) \to \pi_{B}(G) .\] 
Then there exists a map
$\phi^{\prime\prime} \colon Z \to \pi_{B}^{-1}\left(\pi_{B}(G)\right)$
such that $\phi^{\prime\prime}|Z_{0} = \varphi_{0}$ and
$\pi_{B}\phi^{\prime\prime} = \varphi$. Since $f$ is
$0$-invertible (and $\dim Z = 0$), there exists a map
$\phi^{\prime} \colon Z \to Y$ such that
$f\phi^{\prime} = \phi^{\prime\prime}$. Now let
$\phi = g\phi^{\prime}$. Since $g|f^{-1}(G) = f|f^{-1}(G)$,
we have $\varphi_{0} = \phi |Z_{0}$. Finally observe that
the admissibility of $B$ implies $\varphi = \pi_{B}\phi$ as required.

{\em Claim 3. For each countable subset $C \subseteq A$
there exists a countable admissible subset $B \subseteq A$
such that $C \subseteq B$}.

Since $w(Y) = \tau$ and $\dim Y = 0$, it follows
(consult \cite[Theorem 1.3.10]{chibook96}) that $Y$
can be represented as the limit space of a factorizing
$\omega$-spectrum ${\mathcal S}_{Y} = \{ Y_{B}, q^{B}_{C},
\exp_{\omega}A \}$ consisting of zero-dimensional Polish
spaces $Y_{B}$, $B \in \exp_{\omega}A$, and continuous
surjections $q_{C}^{B} \colon Y_{B} \to Y_{C}$, $C \subseteq B$,
$C,B \in \exp_{\omega}A$. Consider also the standard factorizing
$\omega$-spectrum $\displaystyle {\mathcal S}_{X} =
\{ \prod\{ X_{a} \colon a \in B\} ,\pi_{C}^{B}, \exp_{\omega}A \}$
consisting of countable subproducts of the product
$\displaystyle \prod\{ X_{a} \colon a \in A \}$ and corresponding
natural projections. Obviously the full product coincides
with the limit of ${\mathcal S}_{X}$. One more factorizing
$\omega$-spectrum arises naturally. This is the spectrum
${\mathcal S}_{G} = \{ \pi_{B}(G), \pi_{C}^{B}|\pi_{B}(G),
\exp_{\omega}A \}$ the limit of which coincides with $G$. 

Consider the map $f \colon \lim {\mathcal S}_{Y} \to
\lim {\mathcal S}_{X}$. By \cite[Theorem 1.3.4]{chibook96},
there is a cofinal and $\omega$-closed subset
${\mathcal T}_{f}$ of $\exp_{\omega}A$ such that for each
$B \in {\mathcal T}_{f}$ there is a map
$f_{B} \colon Y_{B} \to \prod\{ X_{a} \colon a \in B \}$ such
that $f_{B}\circ q_{B} = \pi_{B}\circ f$. Moreover, these maps
form a morphism
\[ \{ f_{B} ; B \in {\mathcal T}_{f}\} \colon {\mathcal S}_{Y}
\to {\mathcal S}_{X}\]
limit of which coincides with $f$. Since $f$ is proper and
functionally closed, we may assume (see
\cite[Proposition 6.2.9]{chibook96}) without loss of generality
(considering a smaller cofinal and $\omega$-subset of
${\mathcal T}_{f}$ if necessary) that the above indicated morphism
is bicommutative. This simply means that
$q_{B}f^{-1}(K) = f_{B}^{-1}\left(\pi_{B}(K)\right)$ for any
$B \in {\mathcal T}_{f}$ and any closed subset $K$ of the product
$\displaystyle \prod\{ X_{a} \colon a \in A\}$.

Similarly, applying \cite[Theorem 1.3.4]{chibook96} to the
map $g \colon \lim {\mathcal S}_{Y} \to \lim {\mathcal S}_{G}$,
we obtain a cofinal and $\omega$-closed subset
${\mathcal T}_{g}$ of $\exp_{\omega}A$ and the associated to
it morphism
\[ \{ g_{B} \colon Y_{B} \to \pi_{B}(G) ; B \in
{\mathcal T}_{g} \} \colon {\mathcal S}_{Y} \to {\mathcal S}_{G}\]

\noindent limit of which coincides with the map $g$.

By Proposition \ref{P:3.1.1}, the intersection
${\mathcal T} = {\mathcal T}_{f} \cap {\mathcal T}_{g}$
is still a cofinal and $\omega$-closed subset of $\exp_{\omega}A$.
It therefore suffices to show that each $B \in {\mathcal T}$ is
an admissible subset of $A$. Consider a point
$x \in \pi_{B}^{-1}\left( \pi_{B}(G)\right)$. First observe
that the bicommutativity of the morphism associated with
${\mathcal T}_{f}$ implies that
$q_{B}(f^{-1}(x)) = f_{B}^{-1}(\pi_{B}(x))$. Since the
maps $f_{B}$ and $g_{B}$ coincide on $f_{B}^{-1}(\pi_{B}(G))$ we have 
\begin{multline*}
\pi_{B}(g(f^{-1}(x))) = g_{B}(q_{B}(f^{-1}(x))) =
g_{B}(f_{B}^{-1}(\pi_{B}(x))) =\\ f_{B}(f_{B}^{-1}(\pi_{B}(x))) =
\pi_{B}(x)
\end{multline*}
as required.

{\em Claim 4. If $C$ and $B$ are admissible subsets of $A$
and $C \subseteq B$, then the map
$\pi_{C}^{B}|\pi_{B}(G) \colon \pi_{B}(G) \to \pi_{C}(G)$
is $0$-soft}.

This property follows from Claim 2 and \cite[Lemma 6.1.15]{chibook96}.

After having all the needed properties of admissible subsets
established we proceed as follows. Since $|A| = \tau$ we can
write $A = \{ a_{\alpha} \colon \alpha < \tau\}$. By Claim 3,
each $a_{\alpha} \in A$ is contained in a countable admissible
subset $B_{\alpha} \subseteq A$. Let
$A_{\alpha} = \cup\{ B_{\beta} \colon \beta \leq \alpha\}$.
We use these sets to define a transfinite inverse spectrum
${\mathcal S} = \{ G_{\alpha}, p_{\alpha}^{\alpha +1}, \tau \}$
as follows. Let $G_{\alpha} = \pi_{A_{\alpha}}(G)$ and
$p_{\alpha}^{\alpha +1} = \pi_{A_{\alpha}}^{A_{\alpha +1}}|
G_{\alpha +1}$ for each $\alpha < \tau$. All the required
properties of the spectrum ${\mathcal S}_{G}$ are satisfied by
construction. 

The implication (b) $\Longrightarrow$ (a) immediately
follows from \cite[Proposition 6.3.4]{chibook96}
\end{proof}

\begin{rem}
Actually $0$-soft homomorphism $p_{\alpha}^{\alpha +1} \colon
G_{\alpha+1} \to G_{\alpha}$, $\alpha < \tau$, in Theorem \ref{T:ae0}(b)1
has Polish kernel in a somewhat
 stonger sense than the original definition  presented in
Subsection \ref{SS:aen}. Namely a Polish space $P$ (from
the definition), such that $G_{\alpha +1}$ admist a
$C$-embedding into the product $G_{\alpha} \times P$ in such a way that
$p_{\alpha}^{\alpha +1}$ coincides with the restriction of
the projection
$\pi_{G_{\alpha}} \colon G_{\alpha} \times P \to G_{\alpha}$,
can be chosen to be a Polish group and the embedding of
$G_{\alpha +1} \to G_{\alpha} \times P$ can be assumed to
be a homomorphism. Of course this implies that
$\ker p_{\alpha}^{\alpha +1}$, as a closed subgroup of $P$, is
itself a Polish group. It would be interesting to see whether
the converse of this observation is also true, i.e. is it true that
if the kernel $\ker p$ of a $0$-soft homomorphism $p \colon G \to L$
of $\text{AE}(0)$-groups is Polish, then $p$ has a Polish kernel
in the sense of Subsection \ref{SS:aen}.
\end{rem}

Next we characterize $0$-soft homomorphisms of $\text{AE}(0)$-groups
with Polish kernels.

\begin{pro}\label{P:polishkernel}
A $0$-soft homomorphism $f \colon G \to L$ between $\text{AE}(0)$-groups
has a Polish kernel
if and only if there exist factorizing $\omega$-spectra\;
${\mathcal S}_{G} = \{ G_{\alpha}, p_{\alpha}^{\beta}, A \}$,
${\mathcal S}_{L} = \{ L_{\alpha}, q_{\alpha}^{\beta}, A \}$,
consisting of Polish groups and $0$-soft limit homomorphisms, and
a morphism 
$\{ f_{\alpha} \} \colon {\mathcal S}_{G} \to {\mathcal S}_{L}$,
consisting of $0$-soft homomorphisms, such that the following
conditions are satisfied:
\begin{enumerate}
\item
$G = \lim{\mathcal S}_{G}$, $L = \lim{\mathcal S}_{L}$ and
$f = \lim \{ f_{\alpha} \}$.
\item
All limit projections of the spectra ${\mathcal S}_{G}$ and
${\mathcal S}_{L}$ are $0$-soft.
\item
All limit square diagrams, generated by limit projections of
spectra ${\mathcal S}_{G}$ and ${\mathcal S}_{L}$, by elements
of the morphism $\{ f_{\alpha} \}$ and by the map $f$, are the
Cartesian squares.
\end{enumerate}
\end{pro}
\begin{proof}
By Proposition \ref{P:B8.2.1}, we may assume, without loss of
generality, that
both $G$ and $L$ are topologically and algebraically isomorphic
to the limits of factorizing $\omega$-spectra
${\mathcal S}_{G} = \{ G_{\alpha}, p_{\alpha}^{\beta}, A \}$
and ${\mathcal S}_{L} = \{ L_{\alpha}, q_{\alpha}^{\beta}, A \}$,
consisting of Polish groups and $0$-soft limit homomorphisms.
By \cite[Theorem 1.3.6]{chibook96}, we may also assume that the
map
$f$ is the limit of a morphism
$\{ f_{\alpha} \colon G_{\alpha} \to L_{\alpha} ; A\}$,
consisting of continuous maps
$f_{\alpha} \colon G_{\alpha} \to L_{\alpha}$.
Since $f$ itself and all the limit projections
$p_{\alpha} \colon G \to G_{\alpha}$ and
$q_{\alpha} \colon L \to L_{\alpha}$ are homomorphisms
between respective groups, it follows easily that
$f_{\alpha} \colon G_{\alpha} \to L_{\alpha}$ is
also a homomorphism, $\alpha \in A$. Finally since $f$ is
$0$-soft and has a Polish kernel, it follows, by 
\cite[Theorem 6.3.1(vi)]{chibook96}, that all limit square diagrams

\[
\begin{CD}
G @>f>> L\\
@V{p_{\alpha}}VV @VV{q_{\alpha}}V\\
G_{\alpha} @>f_{\alpha}>> L_{\alpha},
\end{CD}
\]
\bigskip

\noindent generated by limit projections of
spectra ${\mathcal S}_{G}$ and ${\mathcal S}_{L}$, by elements
of the morphism $\{ f_{\alpha} \}$ and by the map $f$, are the
Cartesian squares.
\end{proof}

Below, in Subsection \ref{SS:universal}, we consider actions of
$\text{AE}(0)$-groups on $\text{AE}(0)$-spaces. Main tool here
is the following statement (see \cite[Theorem 8.7.1]{chibook96}).

\begin{pro}\label{P:action}
Let $\lambda \colon G \times X \to X$ be a continuous action of
an $\text{AE}(0)$-group $G$ 
on an $\text{AE}(0)$-space $X$. Suppose that $X$ is homeomorphic
to the limit space of a factorizing $\omega$-spectrum
${\mathcal S}_{X} = \{ X_{\alpha}, p_{\alpha}^{\beta}, A\}$
consisting of Polish spaces and $0$-soft limit projections.
Suppose also that $G$ is topologically and algebraically
isomorphic to the limit of the factorizing $\omega$-spectrum
${\mathcal S}_{G} = \{ G_{\alpha}, s_{\alpha}^{\beta}, A\}$
consisting of Polish groups and $0$-soft limit homomorphisms.
Then $\lambda$ is the limit of ``level actions", i.e.
$\lambda = \lim \lambda_{\alpha}$, where 
\[ \{ \lambda_{\alpha} \colon G_{\alpha} \times X_{\alpha}
\to X_{\alpha} , B\} \colon {\mathcal S}_{G}|B \times
{\mathcal S}_{X}|B \to {\mathcal S}_{X}|B\]
 \noindent is a morphism between the spectra
${\mathcal S}_{G}|B \times {\mathcal S}_{X}|B$ and
${\mathcal S}_{X}|B$ and $B$ is a cofinal and
$\omega$-closed subset of the indexing set $A$.
\end{pro}

\section{Applications}\label{S:applications}

\subsection{Universal $\text{AE}(0)$-groups and universal
actions of $\text{AE}(0)$-groups}\label{SS:universal}
In this Subsection we prove the existence of universal
$\text{AE}(0)$-groups of a given weight as well as the
existence of a universal action of a $\text{AE}(0)$-group
of a given weight on a compact $\text{AE}(0)$-space of the
same weight.

\begin{pro}\label{P:univ}
Let $\tau \geq \omega$. The class of $\text{AE}(0)$-groups of weight
$\leq \tau$ contains a universal element. More formally, every
$\text{AE}(0)$-group is topologically and algebraically isomorphic
to a closed and $\text{C}$-embedded subgroup of the power
$\left(\operatorname{Aut}({\mathbb Q})\right)^{\tau}$,
where $\operatorname{Aut}({\mathbb Q})$
denotes the group of autohomeomorphisms of the Hilbert cube ${\mathbb Q}$.
\end{pro}
\begin{proof}
Let $G$ be a $\text{AE}(0)$-group of weight $\tau$. By
Corollary \ref{C:emb}, $G$ is topologically and algebraically
isomorphic to a closed and $C$-embedded subgroup of the
product $\displaystyle \prod\{ G_{t} \colon t \in T\}$,
where $G_{t}$ is a Polish group for each $t \in T$ and $|T| = \tau$.
By Uspenskii's theorem \cite{usp86}, $G_{t}$ can be identified with a closed 
subgroup of $\operatorname{Aut}({\mathbb Q})$. Obviously
$\displaystyle \prod\{ G_{t} \colon t \in T\}$,
and consequently $G$, is a closed and
$C$-embedded subgroup of
$\left(\operatorname{Aut}({\mathbb Q})\right)^{\tau}$.
\end{proof}

Let $\tau > \omega$. Clearly the $\text{AE}(0)$-group
$\left(\operatorname{Aut}({\mathbb Q})\right)^{\tau}$ (i.e. the $\tau$-th power
of the group $\operatorname{Aut}({\mathbb Q})$)
continuously acts on the Tychonov cube ${\mathbb Q}^{\tau}$
via the natural action (``coordinatewise evaluation")
\[ \operatorname{ev}_{\tau} \colon \left(\operatorname{Aut}
({\mathbb Q})\right)^{\tau} \times
{\mathbb Q}^{\tau} \to {\mathbb Q}^{\tau} ,\]

\noindent which is defined by letting 
\begin{multline*}
\operatorname{ev}_{\tau}\left( \{ g_{\alpha} \colon \alpha < \tau \} ,
\{ q_{\alpha} \colon \alpha < \tau\}\right) =
\left\{ g_{\alpha}(q_{\alpha}) \colon \alpha < \tau \right\}
\;\;\text{for each}\\
 \left( \{ g_{\alpha} \colon \alpha < \tau \} ,
\{ q_{\alpha} \colon \alpha < \tau\}\right) \in
\left(\operatorname{Aut}({\mathbb Q})\right)^{\tau} \times
{\mathbb Q}^{\tau} .
\end{multline*}

\begin{thm}\label{T:actionuniv}
Let $\tau > \omega$. The action $\operatorname{ev}_{\tau} \colon
\left(\operatorname{Aut}({\mathbb Q})\right)^{\tau} \times
{\mathbb Q}^{\tau} \to {\mathbb Q}^{\tau}$ is universal
in the category of actions of $\text{AE}(0)$-groups
of weight $\tau$ on  compact $\text{AE}(0)$-spaces of weight
$\tau$. More formally,
let $\lambda \colon G \times X \to X$ be a continuous action of a
$\text{AE}(0)$-group $G$ of weight $\tau$ on a compact
$\text{AE}(0)$-space $X$ of weight
$\tau$. Then there exists a topological and algebraic embedding
$i_{G} \colon G \to \left(\operatorname{Aut}
({\mathbb Q})\right)^{\tau}$ with a closed image and an embedding
$i_{X} \colon X \to {\mathbb Q}^{\tau}$ such that the following diagram
\[
\begin{CD}
\left(\operatorname{Aut}
({\mathbb Q})\right)^{\tau}\times {\mathbb Q}^{\tau}
@>\operatorname{ev}_{\tau}>> {\mathbb Q}^{\tau}\\
@A{i_{G}}\times i_{X}AA @AA{i_{X}}A\\
G \times X @>\lambda>> X
\end{CD}
\]

\noindent commutes.
\end{thm}
\begin{proof}
By Proposition \ref{P:B8.2.1},  $G$ can be represented as
the limit of a factorizing $\omega$-spectrum
${\mathcal S}_{G} = \{ G_{\alpha}, s_{\alpha}^{\beta}, A\}$
consisting of Polish groups and $0$-soft limit homomorphisms.
Similarly, by \cite[Proposition 6.3.5]{chibook96}, $X$ can be
represented as the limit space of a factorizing $\omega$-spectrum
${\mathcal S}_{X} = \{ X_{\alpha}, p_{\alpha}^{\beta}, A\}$
consisting of metrizable compacta and $0$-soft limit projections.
Without loss of generality we may assume that these spectra
${\mathcal S}_{G}$ and ${\mathcal S}_{X}$ have the same indexing
set $A$ and $|A| = \tau$. By Proposition \ref{P:action}, the given
action $\lambda \colon G \times X \to X$ is the limit of level
actions, i.e.
$\lambda = \lim \lambda_{\alpha}$, where 
\[ \{ \lambda_{\alpha} \colon G_{\alpha} \times X_{\alpha}
\to X_{\alpha} , B\} \colon {\mathcal S}_{G}|B \times
{\mathcal S}_{X}|B \to {\mathcal S}_{X}|B\]
 \noindent is a morphism between the spectra
${\mathcal S}_{G}|B \times {\mathcal S}_{X}|B$ and
${\mathcal S}_{X}|B$ and $B$ is a cofinal and
$\omega$-closed subset of the indexing set $A$. We may also
assume that $|B| = \tau$. 

Since $G = \lim {\mathcal S}_{G}|B$ it follows that the
diagonal product

\[ s = \triangle \{ s_{\alpha} \colon G \to G_{\alpha} ,
\alpha \in B\} \colon G \to \prod\{ G_{\alpha} \colon \alpha \in B\}\]

\noindent is a topological and algebraic isomorphism with
a closed image. Similarly
the diagonal product

\[ p = \triangle \{ p_{\alpha} \colon X \to X_{\alpha} ,
\alpha \in B\} \colon X \to \prod\{ X_{\alpha} \colon \alpha \in B\}\]

\noindent is an embedding. Consider also the product action
\[ \widetilde{\lambda} \colon \prod\{ G_{\alpha} \colon \alpha \in B\}
\times \prod\{ X_{\alpha} \colon \alpha \in B\} \to
\prod\{ X_{\alpha} \colon \alpha \in B\}\]

\noindent defined by letting

\begin{multline*}
 \widetilde{\lambda}\left( \{ g_{\alpha} \colon
\alpha \in B\} ,\{ x_{\alpha} \colon \alpha \in B \}\right) =
\{ \lambda_{\alpha}(g_{\alpha}, x_{\alpha}) \colon \alpha \in B\}
\;\;\text{for each}\\
 \left( \{ g_{\alpha} \colon
\alpha \in B\} ,\{ x_{\alpha} \colon \alpha \in B \}\right)
\in \prod\{ G_{\alpha} \colon \alpha \in B\} \times
\prod\{ X_{\alpha} \colon \alpha \in B\} .
\end{multline*}

\noindent Note that
$\widetilde{\lambda}\circ (s\times p) = p \circ \lambda$,
i.e. the following diagram

\[
\begin{CD}
\prod\{ G_{\alpha} \colon \alpha \in B\} \times
\prod\{ X_{\alpha} \colon \alpha \in B\} @>\widetilde{\lambda}>>
\prod\{ X_{\alpha} \colon \alpha \in B\}\\
@A{s\times p}AA @AA{p}A\\
G\times X @>\lambda>> X
\end{CD}
\]

\noindent is commutative.

Since for each $\alpha \in B$ the group $G_{\alpha}$ is Polish, it follows
by \cite{megr96} (see also \cite[Theorem 2.6.7]{bekech96}) that
there exist a
topological and algebraic embedding $j_{\alpha} \colon G_{\alpha} \to \operatorname{Aut}({\mathbb Q}_{\alpha})$
with a closed image and an embedding
$i_{\alpha} \colon X_{\alpha} \to {\mathbb Q}_{\alpha}$
(here ${\mathbb Q}_{\alpha}$ denotes a copy of the
Hilbert cube ${\mathbb Q}$) such that
$\operatorname{ev}_{\alpha} \circ (j_{\alpha}
\times i_{\alpha}) = i_{\alpha} \circ \lambda_{\alpha}$.
This simply means that the diagram

\[
\begin{CD}
\operatorname{Aut}({\mathbb Q}_{\alpha}) \times
{\mathbb Q}_{\alpha} @>\operatorname{ev}_{\alpha}>>
{\mathbb Q}_{\alpha}\\
@A{j_{\alpha}\times i_{\alpha}}AA @AA{i_{\alpha}}A\\
G_{\alpha}\times X_{\alpha} @>\lambda_{\alpha}>> X_{\alpha}
\end{CD}
\]
\medskip

\noindent commutes for each $\alpha \in B$. Here $\operatorname{ev}_{\alpha}
\colon \operatorname{Aut}({\mathbb Q}_{\alpha}) \times
{\mathbb Q}_{\alpha} \to {\mathbb Q}_{\alpha}$ is the evaluation action, i.e.
$\operatorname{ev}_{\alpha}(g_{\alpha}, x_{\alpha}) = g_{\alpha}(x_{\alpha})$ for each $(g_{\alpha}, x_{\alpha}) \in \operatorname{Aut}({\mathbb Q}_{\alpha}) \times
{\mathbb Q}_{\alpha}$. Let 
\[ j  = \times \{ j_{\alpha} \colon \alpha \in B \} \colon
\prod\{ G_{\alpha} \colon \alpha \in B\} \to
\prod\{ \operatorname{Aut}({\mathbb Q}_{\alpha}) \colon \alpha \in B\} \]

\noindent and

\[ i = \times\{ i_{\alpha} \colon \alpha \in B\} \colon
\prod\{ X_{\alpha} \colon \alpha \in B\} \to \prod\{
{\mathbb Q}_{\alpha} \colon \alpha \in B\} .\]

Finally consider the commutative diagram

\[
\begin{CD}
\prod\{ \operatorname{Aut}({\mathbb Q}_{\alpha}) \colon
\alpha \in B\} \times \prod\{
{\mathbb Q}_{\alpha} \colon \alpha \in B\} @>\times\{
\operatorname{ev}_{\alpha} \colon \alpha \in B\}>> \prod\{
{\mathbb Q}_{\alpha} \colon \alpha \in B\}\\
@A{j \times i}AA @AA{i}A\\
\prod\{ G_{\alpha} \colon \alpha \in B\} \times
\prod\{ X_{\alpha} \colon \alpha \in B\} @>\widetilde{\lambda}>>
\prod\{ X_{\alpha} \colon \alpha \in B\}\\
@A{s\times p}AA @AA{p}A\\
G\times X @>\lambda>> X
\end{CD}
\]

\noindent and note that since $|B| = \tau$ the upper horizontal
arrow is actually the action
$\operatorname{ev}_{\tau} \colon \left(\operatorname{Aut}(
{\mathbb Q})\right)^{\tau}\times {\mathbb Q}^{\tau}
\to {\mathbb Q}^{\tau}$. Clearly it suffices to let
$i_{G} = (j \times i) \circ (s \times p)$ and $i_{X} = i \circ p$.
This completes the proof.
\end{proof}

It would be very interesting to prove that for a $\text{AE}(0)$-group
$G$ of weight $\tau > \omega$ the
category of  $AE(0)$-spaces (of weight $\tau$) admitting actions of the
group $G$ and their $G$-maps contains a universal object. For
$\tau = \omega$ this fact has recently been proved in \cite{hjo99}.

\subsection{Closed subgroups of powers of the
symmetric group $S_{\infty}$}\label{SS:symmetric}

The following result gives an embeddability criterion into
the symmetric group $S_{\infty}$ - the group
of all bijections of ${\mathbb N}$ under the relative topology
inherited from ${\mathbb N}^{\mathbb N}$. It is important to note
\cite{dough94} that there exist zero-dimensional Polish groups
which can not be embedded into $S_{\infty}$ as closed subgroups.

\begin{thm}[\cite{bekech96}]\label{T:beke}
Let $G$ be a Polish group. Then the following conditions are equivalent:
\begin{itemize}
\item[(i)]
$G$ is isomorphic to a closed subgroup of $S_{\infty}$;
\item[(ii)]
$G$ admits a (countable) neighborhood basis at the identity
consisting of open subgroups;
\item[(iii)]
$G$ admits a (countable) basis closed under left multiplication
(or a countable basis closed under right multiplication);
\item[(iv)]
$G$ admits a compatible left-invariant ultrametric.
\end{itemize}
\end{thm}

Next we characterize those topological $AE(0)$-groups
which are isomorphic to closed subgroups of infinite powers
$S_{\infty}^{\tau}$, $\tau \geq \omega$, of $S_{\infty}$.

\begin{thm}\label{T:closed}
Let $\tau \geq 1$ be a cardinal number. The following conditions are equivalent for any
topological $AE(0)$-group $G$ of weight $\tau \geq \omega$:
\begin{itemize}
\item[(i)]
$G$ is isomorphic to a closed subgroup of $S_{\infty}^{\tau}$;
\item[(ii)]
$G$ admits a neighborhood basis at the identity
consisting of open subgroups.
\end{itemize}
\end{thm}
\begin{proof}
If $\tau = 1$ our statement coincides with Theorem \ref{T:beke}. Next
consider the case $1 < \tau \leq\omega$. It is easy to see that the group
$S_{\omega}^{\tau}$ admits a countable neighborhood basis
at the identity consisting of open subgroups. Obviously 
every closed subgroup of $S^{\omega}_{\tau}$ has the same
property (and, consequently, by Theorem \ref{T:beke}, can be
embedded into $S_{\infty}$). Conversely if a Polish group $G$ admits a countable
neighborhood basis at the identity
consisting of open subgroups, then, by Theorem \ref{T:beke},
$G$ is isomorphic to a closed subgroup of $S_{\infty}$. It only
remains to note that $G_{\infty}$ is isomorphic to a closed
subgroup of $S_{\infty}^{\tau}$ for any $\tau$.

Next we assume that $\tau > \omega$. By
\cite[Lemma 8.2.1]{chibook96}, $G$ is isomorphic
to the limit space of a factorizing $\omega$-spectrum
${\mathcal S}_{G} = \{ G_{\alpha}, p_{\alpha}^{\beta}, A\}$
all spaces $G_{\alpha}$, $\alpha \in A$, of which are Polish
groups and all limit projections $p_{\alpha} \colon G \to G_{\alpha}$,
$\alpha \in A$,
of which are $0$-soft homomorphisms.

Now consider the following relation $L \in A^{2}$:
\begin{multline*}
 L = \left\{ (\alpha ,\beta ) \in A^{2} \colon \alpha \leq \beta
\;\text{and there exists a countable neighborhood basis}\right.\\
 \;\;\;\;\;\;\;\; {\mathcal V}_{\alpha} \;\text{at
the identity}\; e_{\alpha} \;\text{of}\; G_{\alpha} \; \text{containing
intersections of its finite subcoll-}\\
\;\;\;\;\;\;\;\;\text{ections and such that for each}\; V \in {\mathcal V}_{\alpha}\;\;
\text{there is an open subgroup} \; U^{\beta ,\alpha}_{V}\\ 
\left.  \; \text{of}\; G_{\beta}\;
\text{with}\; U_{V}^{\beta ,\alpha} \subseteq
\left( p_{\alpha}^{\beta}\right)^{-1}(V)\; \right\}\;\;\;\;\;\;\;\;\;\;\;\;\;\;\;\;\;\;\;\;\;\;\;\;\;\;
\;\;\;\;\;\;\;\;\;\;\;\;\;\;\;\;\;\;\;\;\;\;\;\;\;
\;\;\;\;\;\;\;
\end{multline*}

Let us verify conditions of Proposition \ref{P:search}.

{\em Existence}.
For each $\alpha \in A$ we need to find
$\beta \in A$  such that  $(\alpha ,\beta ) \in L$.
Let ${\mathcal V}_{\alpha}$ be a countable
neighborhood basis at $e_{\alpha} \in G_{\alpha}$ which contains
intersections of its finite subcollections. For each $V \in {\mathcal V}_{\alpha}$
the set $p_{\alpha}^{-1}\left( V\right)$ is a
neighborhood of the identity $e \in G$. By (ii), there exists an
open subgroup $U_{V}$ of $G$ such that
$U_{V} \subseteq p_{\alpha}^{-1}\left( V\right)$.
Since every open subgroup in $G$ is closed, it follows that $U_{V}$, as
an open and closed subset of $G$, is a functionally open in $G$.
Recall that the spectrum
${\mathcal S}_{G}$
is factorizing and consequently there exist an index
$\beta_{V} \in A$ and an
open subset $U_{V}^{\beta_{V}} \subseteq G_{\beta_{V}}$ such
that $\beta_{V} \geq \alpha$
and $U_{V} = p_{\beta_{V}}^{-1}\left( U_{V}^{\beta_{V}}\right)$.
By Corollary \ref{C:3.1.2}, there exists an index $\beta \in A$ such that
$\beta \geq \beta_{V}$ for each $V \in {\mathcal V}_{\alpha}$.
Let $U_{V}^{\beta ,\alpha} = \left( p_{\beta_{V}}^{\beta}\right)^{-1}\left( U_{V}^{\beta_{V}}\right)$,
$V \in {\mathcal V}_{\alpha}$. Note that

\begin{multline*}
U_{V}^{\beta ,\alpha} =
\left( p_{\beta_{V}}^{\beta}\right)^{-1}\left( U_{V}^{\beta_{V}}\right)
= p_{\beta}\left( p_{\beta_{V}}\left( U_{V}^{\beta_{V}}\right)\right)
= p_{\beta}(U_{V})
\subseteq p_{\beta}\left( p_{\alpha}^{-1}(V)\right) =\\
 \left( p_{\alpha}^{\beta}\right)^{-1}(V)
\;\text{for each}\; V \in {\mathcal V}_{\alpha} .
\end{multline*}

\noindent Since the limit projection $p_{\beta} \colon G \to G_{\beta}$ is a
homomorphism it follows that $U_{V}^{\beta ,\alpha} = p_{\beta}(U_{V})$
is a subgroup of $G_{\beta}$. Finally since $U_{V}^{\beta_{V}}$
is open in $G_{\beta_{V}}$, we conclude that $U_{V}^{\beta ,\alpha} = \left( p_{\beta_{V}}^{\beta}\right)^{-1}\left( U_{V}^{\beta_{n}}\right)$
is open in $G_{\beta}$. This shows that $(\alpha ,\beta ) \in L$.

{\em Majorantness}.
If  $(\alpha ,\beta ) \in L$  and
$\gamma \geq \beta$, then  $(\alpha ,\gamma ) \in L$.
Since $(\alpha ,\beta ) \in L$ it follows that for each
$V \in {\mathcal V}_{\alpha}$ there exists an open subgroup $U_{V}^{\beta ,\alpha}$
of $G_{\beta}$ such that
$U_{V}^{\beta ,\alpha} \subseteq \left( p_{\alpha}^{\beta}\right)^{-1}
(V)$ where
${\mathcal V}_{\alpha}$ is a countable
neighborhood basis at the identity $e_{\alpha} \in G_{\alpha}$
containing intersections of its finite subcollections.
Let $U_{V}^{\gamma ,\alpha} =
\left( p_{\beta}^{\gamma}\right)^{-1}(U_{V}^{\beta ,\alpha})$ for
each $V \in {\mathcal V}_{\alpha}$. Since the projection
$p_{\beta}^{\gamma} \colon G_{\gamma} \to G_{\beta}$ is a
continuous homomorphism, it follows that $U_{V}^{\gamma ,\alpha}$ is an
open subgroup of $G_{\gamma}$. Obviously 

\[ U_{V}^{\gamma ,\alpha} = \left( p_{\beta}^{\gamma}\right)^{-1}(U_{V}^{\beta ,\alpha})
\subseteq  \left( p_{\beta}^{\gamma}\right)^{-1}\left(\left(
p_{\alpha}^{\beta}\right)^{-1}(V)\right) = 
\left( p_{\alpha}^{\gamma}\right)^{-1}(V) ,\;\; V \in {\mathcal V}_{\alpha} , \]

\noindent which shows that $(\alpha ,\gamma ) \in L$ as required.

{\em $\omega$-closeness}.
Let $\{ \alpha_{i} : i \in \omega \}$
be a countable chain in $A$ and
$(\alpha_{i}, \beta ) \in L$ for some
$\beta \in A$ and each $i \in \omega$. We need to show that
$(\alpha ,\beta ) \in L$ where $\alpha =
\sup \{\alpha_{i} \colon i \in \omega \}$.

Let ${\mathcal V}_{\alpha_{i}}$ be a countable
neighborhood basis 
at $e_{\alpha_{i}} \in G_{\alpha_{i}}$ and
$U_{V}^{\beta ,\alpha_{i}}$, $V \in {\mathcal V}_{\alpha_{i}}$,
be an open subgroup of $G_{\beta}$ witnessing the fact that
$(\alpha_{i}, \beta ) \in L$.

Consider the collection 
\[ {\mathcal V}_{\alpha} =
\bigcup\{ \left( p_{\alpha_{i}}^{\alpha}\right)^{-1}\left(
{\mathcal V}_{\alpha_{i}}\right)
\colon i \in \omega\}  .\]

\noindent Since the spectrum ${\mathcal S}_{G}$ is an $\omega$-spectrum and
since $\alpha = \sup\{ \alpha_{i} \colon i \in \omega\}$ it
follows that ${\mathcal V}_{\alpha}$ forms a neighborhood basis
at $e_{\alpha} \in G_{\alpha}$. For each $\widetilde{V} \in {\mathcal V}_{\alpha}$
choose $V \in {\mathcal V}_{\alpha_{i}}$ such that $\widetilde{V} = \left( p_{\alpha_{i}}^{\alpha}\right)^{-1}(V)$ and let $U_{\widetilde{V}}^{\beta ,\alpha} = U_{V}^{\beta ,\alpha_{i}}$. Since $(\alpha_{i} ,\beta) \in L$, we have

\[ U_{\widetilde{V}}^{\beta ,\alpha} = U_{V}^{\beta ,\alpha_{i}} \subseteq \left( p_{\alpha_{i}}^{\beta}\right)^{-1}(V) = \left( p_{\alpha}^{\beta}\right)^{-1}
\left( \left( p_{\alpha_{i}}^{\alpha}\right)^{-1}(V)\right) = \left( p_{\alpha}^{\beta}\right)^{-1}(\widetilde{V}) .\]

\noindent This proves that $(\alpha ,\beta ) \in L$.

According to Proposition \ref{P:search} the set $\widetilde{A}$ of
$L$-reflexive elements is cofinal and $\omega$-closed in $A$. The
$L$-reflexivity of an element $\alpha \in A$ means precisely
that there exists a
countable neighborhood basis ${\mathcal V}_{\alpha}$ at
$e_{\alpha} \in G_{\alpha}$ containing intersections of its finite
subcollections and such that for each $V \in {\mathcal V}_{\alpha}$
there exists an open subgroup $U_{V}^{\alpha} \subseteq G_{\alpha}$
with $U_{V}^{\alpha} \subseteq V$. This obviously means that for each
$\alpha \in {\widetilde{A}}$ the Polish group $G_{\alpha}$
satisfies condition (ii)
of Theorem \ref{T:beke}. Consequently, by Theorem \ref{T:beke},
$G_{\alpha}$ is topologically isomorphic to a closed subgroup
of $S_{\infty}$. 
Next note that since $\widetilde{A}$ is cofinal in $A$ the limit
space of the spectrum
${\mathcal S} = \{ G_{\alpha}, p_{\alpha}^{\beta}, \widetilde{A}\}$
is topologically isomorphic to $G$. This obviously implies that $G$
is isomorphic to a closed subgroup of the product
$\displaystyle \prod\{ G_{\alpha} \colon \alpha \in \widetilde{A}\}$,
which in turn is topologically isomorphic to a closed
subgroup of $S_{\infty}^{\tau}$ (note that
$|\widetilde{A}| = w(G) = \tau$). This completes
the proof of implication (ii) $\Longrightarrow$ (i).

Verification of the implication (i) $\Longrightarrow$ (ii)
is trivial. Proof is completed.
\end{proof}

\begin{cor}\label{C:polish}
Let $\tau \geq 2$. The following conditions are equivalent
for any Polish group $G$:
\begin{itemize}
\item[(i)]
$G$ is topologically isomorphic to a closed subgroup of $S_{\infty}$;
\item[(ii)]
$G$ is topologically isomorphic to a closed subgroup of $S_{\infty}^{\tau}$.
\end{itemize}
\end{cor}

\begin{cor}\label{C:doug}
There exists a zero-dimensional Polish group which is not
topologically isomorphic to a closed subgroup of
$S_{\infty}^{\tau}$ for any cardinal number $\tau$.
\end{cor}
\begin{proof}
It is known \cite{dough94} that there exists a zero-dimensional
Polish group $G$ which is not topologically isomorphic to a
closed subgroup of $S_{\infty}$. By Corollary \ref{C:doug}, $G$
can not be topologically isomorphic to a closed subgroup of
$S_{\infty}^{\tau}$ for
any cardinal $\tau \geq 2$.
\end{proof}


\subsection{Baire isomorphisms}\label{SS:baire}
The main result of this Subsection (Theorem \ref{T:baire}) 
allows us to reduce in many instances (descriptive) set
theoretical considerations of general $\text{AE}(0)$-groups
to those for Polish groups. 

\begin{lem}\label{L:dixmier}
Let $f \colon X \to Y$ be a $0$-soft homomorphism between
$\text{AE}(0)$-groups. Then there
exists a Baire isomorphism $h \colon Y \times \ker f \to X$ such
that $f \circ h = \pi_{Y}$, where
$\pi_{Y} \colon Y \times \ker f \to Y$ stands for the projection
onto the first coordinate.
\end{lem}
\begin{proof}
By Proposition \ref{P:bb}, there exists a Baire map
$g \colon Y \to X$ such that
$f \circ g = \operatorname{id}_{Y}$.
The required Baire isomorphism $h \colon Y \times \ker f \to X$
(not a homomorphism unless $g$
is a homomorphism itself) can now be defined by letting
\[ h(y,a) = g(y)\cdot a ,\;\;\text{for each}\;\; (y,a) \in Y \times \ker f ,\] 

\noindent where $\cdot$ denotes the multiplication operation in $X$. 
\end{proof}

\begin{thm}\label{T:baire}
Every $\text{AE}(0)$-group is Baire isomorphic to the product of Polish groups.
\end{thm}
\begin{proof}
Let $X$ be a $\text{AE}(0)$-group. If $w(X) = \omega$, then $X$ itself is Polish
and there is nothing to prove.

Let now $w(X) = \tau > \omega$. According to Theorem \ref{T:ae0},
$X$ is topologically and algebraically isomorphic to the limit of a
well-ordered continuous spectrum
${\mathcal S}_{X} = \{ X_{\alpha}, p_{\alpha}^{\alpha +1}, \tau \}$
such that $X_{0}$ is a Polish group and the $0$-soft homomorphism
$p_{\alpha}^{\alpha +1} \colon X_{\alpha +1} \to X_{\alpha}$ has a
Polish kernel for each $\alpha < \tau$.

Our goal is to prove that $X$ is Baire isomorphic to the product
$\displaystyle X_{0} \times \prod\{
\ker p_{\alpha}^{\alpha +1} \colon \alpha < \tau\}$. We proceed by induction.
By Lemma \ref{L:dixmier}, there exists a Baire isomorphism
$h_{1} \colon X_{0} \times \ker p_{0}^{1} \to X_{1}$ such
that $p_{0}^{1} \circ h_{1} = \pi_{X_{0}}$.
Suppose that for each $\alpha$, where $1 \leq \alpha < \beta < \tau$, we
have already constructed Baire isomorphism
$\displaystyle h_{\alpha} \colon X_{0} \times \prod\{
\ker p_{\delta}^{\delta +1}\colon \delta < \alpha\} \to X_{\alpha}$
in such a way that
\begin{enumerate}
\item
If $\alpha + 1 < \beta$, then
\[
\begin{CD}
X_{0} \times \prod\{ \ker p_{\delta}^{\delta +1}\colon \delta <
\alpha +1\} @>h_{\alpha +1}>> X_{\alpha +1}\\
@V{\operatorname{id}_{X_{0}}\times \pi_{\alpha}^{\alpha +1}}VV @VV
{p_{\alpha}^{\alpha +1}}V\\
X_{0}\times \prod\{ \ker p_{\delta}^{\delta +1}\colon
\delta < \alpha\} @>h_{\alpha}>> X_{\alpha},
\end{CD}
\]

\noindent where 
\begin{multline*}
 \pi_{\alpha}^{\alpha +1} \colon \prod\{ \ker
p_{\delta}^{\delta +1}\colon \delta < \alpha +1\}  =\\
\prod\{ \ker p_{\delta}^{\delta +1}\colon \delta < \alpha\}
\times X_{\alpha} \to \prod\{ \ker
p_{\delta}^{\delta +1}\colon \delta < \alpha\}
\end{multline*}

\noindent is the natural projection.
\item
$h_{\alpha} = \lim\{ h_{\gamma} \colon \gamma < \alpha\}$, whenever
$\alpha <\beta$ is a limit ordinal number.
\end{enumerate} 

We now construct Baire isomorphism
$\displaystyle h_{\beta} \colon X_{0} \times \prod\{
\ker p_{\delta}^{\delta +1}\colon \delta < \beta\} \to X_{\beta}$

If $\beta$ is a limit ordinal number, then we let
$h_{\beta} = \lim\{ h_{\alpha} \colon \alpha < \beta\}$.

If $\beta = \alpha +1$, then consider the following commutative diagram

\[
\begin{CD}
\left( X_{0} \times \prod\{ \ker p_{\delta}^{\delta +1}\colon \delta <
\alpha\}\right) \times \ker p_{\alpha}^{\alpha +1} @>h_{\alpha}\times \operatorname{id}>> X_{\alpha} \times \ker p_{\alpha}^{\alpha +1} @>h>> X_{\alpha +1}\\
@V{\pi_{1}}VV   @V{\pi_{X_{\alpha}}}VV @VV
{p_{\alpha}^{\alpha +1}}V\\
X_{0}\times \prod\{ \ker p_{\delta}^{\delta +1}\colon
\delta < \alpha\} @>h_{\alpha}>> X_{\alpha}
@>\operatorname{id}_{X_{\alpha}}>> X_{\alpha},
\end{CD}
\]
\bigskip

\noindent where 

\begin{itemize}
\item[(a)]
$\operatorname{id} \colon \ker p_{\alpha}^{\alpha +1} \to
\ker p_{\alpha}^{\alpha +1}$ stands for the identity map;
\item[(b)]
\[ \pi_{1} \colon \left( X_{0} \times \prod\{ \ker
p_{\delta}^{\delta +1}\colon \delta <
\alpha\}\right) \times \ker p_{\alpha}^{\alpha +1}
\to X_{0} \times \prod\{ \ker p_{\delta}^{\delta +1}\colon \delta <
\alpha\}\]

\noindent denotes the projection onto the first coordinate and
\item[(c)]
$h \colon X_{\alpha} \times \ker p_{\alpha}^{\alpha +1} \to
X_{\alpha +1}$ is a Baire isomorphism
existence of which is guaranteed by Lemma \ref{L:dixmier}.
\end{itemize}
The required Baire isomorphism $\displaystyle h_{\alpha +1}
\colon X_{0} \times \prod\{ \ker p_{\delta}^{\delta +1} \colon
\delta < \alpha +1\} \to X_{\alpha +1}$ can now be defined as
the composition $h_{\alpha +1} = h\circ \left( h_{\alpha} \times
\operatorname{id}\right)$.
This completes induction and finishes the construction of Baire
isomorphisms
$h_{\alpha}$, $\alpha < \tau$. It is now easy to see that 
\[ h = \lim\{ h_{\alpha} \colon \alpha < \tau\} \colon
X_{0} \times \prod\{ \ker p_{\alpha}^{\alpha +1} \colon \alpha < \tau \} \to X\]

\noindent is the required Baire isomorphism. Note here that $X_{0}$ as well as
$\ker p_{\alpha}^{\alpha +1}$, $\alpha < \tau$, are Polish groups.
Proof is completed.
\end{proof}

\providecommand{\bysame}{\leavevmode\hbox to3em{\hrulefill}\thinspace}



\end{document}